\documentclass[fleqn]{mat01}
\usepackage{times,mathtimy,amssymb,latexsym}%,rangecite}
\begin{document}

\setcounter{page}{259} \firstpage{259}

\newtheorem{theore}{Theorem}
\renewcommand\thetheore{\arabic{theore}}
\newtheorem{theor}[theore]{\bf Theorem}

\newtheorem{lem}{Lemma}
\newtheorem{case}{Case}

\title{Solutions for a class of iterated singular equations}

\markboth{A \c{C}etinkaya and N \"{O}zalp}{Iterated singular
equations}

\author{A \c{C}ET$\dot{\rm I}$NKAYA and N \"{O}ZALP$^{*}$}

\address{Faculty of Arts and Sciences, Department of Mathematics,
Ahi Evran University, 40100~Kirseh$\dot{\rm i}$r, Turkey\\
\noindent $^{*}$Faculty of Sciences, Department of Mathematics,
Ankara University, Be\c{s}evler, 06100~Ankara, Turkey\\
\noindent E-mail: caysegul@gazi.edu.tr;
nozalp@science.ankara.edu.tr}

\volume{117}

\mon{May}

\parts{2}

\pubyear{2007}

\Date{MS received 15 March 2004; revised 30 March 2007}

\begin{abstract}
Some fundamental solutions of radial type for a class of iterated
elliptic singular equations including the iterated Euler equation
are given.
\end{abstract}


\maketitle

\section{Introduction}

Consider the class of equations
\begin{equation}
Lu=\sum_{i=1}^{n}\left( \frac{r}{x_{i}}\right) ^{p}\left[
x_{i}^{2}\frac{\partial ^{2}u}{\partial x_{i}^{2}}+\alpha
_{i}x_{i}\frac{\partial u}{\partial x_{i}}\right] +\lambda u=0,
\end{equation}
where $\lambda ,\alpha _{i}\,\,(i=1,2,\dots,n)\,$ are real
parameters, $p\,(>0) $ is a real constant and $r$ is defined by
\begin{equation}
r^{p}=x_{1}^{p}+x_{2}^{p}+\dots+x_{n}^{p}.
\end{equation}

The domain of the operator $L$ is the set of all real-valued
functions $u(x)$ of the class $C^{2}(D),$ where
$x=(x_{1},x_{2},\dots,x_{n})$ denotes points in $R^{n}$ and $D$ is
the regularity domain of $u$ in $R^{n}.$ Note that $(1)$ includes
the Laplace equation and an equidimensional (Euler) equation as
special cases.

In \cite{1} and \cite{2}, Alt\i n studied radial type solutions of
a class of singular partial differential equations of even order
and obtained Lord Kelvin principle for this class of equations.

In \cite{4}, all radial type solutions of eq.~(1) are obtained by
showing that for all solutions of the form $u=f(r^{m}),$ $f\in
C^{2},$ the function $f$ satisfies
\begin{equation*}
f(r^{m})=r^{cm},
\end{equation*}%
where $c$ is a root of the equation
\begin{equation*}
m^{2}c^{2}+m\left( -p+n(p-1)+\sum\limits_{i=1}^{n}\alpha
_{i}\right) c+\lambda =0.
\end{equation*}

In \cite{3,5}, \"{O}zalp and \c{C}etinkaya obtained expansion
formulas and Kelvin principle for the iterates of eq.~(1). Lyakhov
and Ryzhkov \cite{6} obtained Almansi's expansions for
\hbox{$B$-polyharmonic} equation i.e. obtained the solutions of
the equation
\begin{equation*}
\Delta _{B}^{m}f=0,
\end{equation*}
where
\begin{equation*}
\Delta _{B}=\sum_{j=1}^{n}B_{j}+\sum_{i=n+1}^{N}\frac{\partial
^{2}}{\partial x_{i}^{2}}, \qquad B_{j}=\partial ^{2}/\partial
x_{j}^{2}+\frac{\gamma _{j}}{x_{j}}\frac{\partial }{\partial
x_{j}}.
\end{equation*}
In this paper, as a continuation of \cite{3}, we consider the
class of equations
\begin{equation}
\left( \prod\limits_{j=1}^{q}L_{j}^{k_{j}}\right) u=(
L_{1}^{k_{1}}L_{2}^{k_{2}}\dots L_{q}^{k_{q}}) u=0,
\end{equation}
where $q,~k_{1},\dots,k_{q}$ are positive integers,
$\lambda_{j},~\alpha _{i}^{(j)}~(j=1,2,\dots,q;~i=1,2,\dots, n)$
are real constants,
\begin{equation*}
L_{j}=\sum_{i=1}^{n}\left( \frac{r}{x_{i}}\right) ^{p}\left[
x_{i}^{2}\frac{\partial ^{2}}{\partial x_{i}^{2}}+\alpha
_{i}^{(j)}x_{i}\frac{\partial }{\partial x_{i}}\right] +\lambda
_{j}
\end{equation*}
and the operator $L_{j}^{k_{j}}$ denotes, as usual, the successive
applications of the operator $L_{j}$ onto itself, that is
$L_{j}^{k_{j}}u=L_{j}(L_{j}^{k_{j}-1}u)$.

\section{Solutions for the iterated equation}

We first give some properties of the operator $L_{j}$ (see
\cite{3,4}). By direct computation, it can be shown that
\begin{equation}
L_{j}(r^{m})=\beta _{j}(m)r^{m},
\end{equation}%
where
\begin{equation}
\beta _{j}(m)=[ m(m+2\phi _{j})+\lambda _{j}]
\end{equation}%
and%
\begin{equation}
2\phi _{j}=-p+n(p-1)+\sum_{i=1}^{n}\alpha _{i}^{(j)}.
\end{equation}

The proof of the following lemma can be done easily by using
induction argument on $k_{j}$. For a special case of the lemma,
see \cite{4}.

\begin{lem}
For any real parameter $m,$
\begin{equation*}
L_{j}^{k_{j}}(r^{m})=\beta _{j}^{k_{j}}(m)r^{m},
\end{equation*}
where the integer $k_{j}$ is the iteration number.
\end{lem}

By the linearity of the operators $L_{j}^{{}}$ and by Lemma~1, we
have the following result.

\begin{lem}
\begin{equation}
\left( \prod\limits_{j=1}^{q}L_{j}^{k_{j}}\right) (r^{m})=\left(
\prod\limits_{j=1}^{q}\beta _{j}^{k_{j}}(m)\right) r^{m}.
\end{equation}
\end{lem}

The following theorem, states a class of solutions for the
iterated equations which is our main result.

\begin{theor}[\!]
The function defined by
\begin{align}
u &= \sum_{v\in I_{1}}\sum_{l=0}^{k_{v}-1}r^{-\phi _{v}}\left[
A_{l}~r^{\sqrt{\phi _{v}^{2}-\lambda _{v}}}+B_{l}~r^{-\sqrt{\phi
_{v}^{2}-\lambda _{v}}}\right] (\ln r)^{l}  \notag \\[.3pc]
&\quad\, +\sum_{v\in I_{2}}\sum_{l=0}^{k_{v}-1}r^{-\phi
_{v}}\!\left[\!C_{l}\cos \left(\!\sqrt{\lambda _{v}\!-\!\phi
_{v}^{2}}\ln r\!\right)\!+\!D_{l}\sin \left(\!\sqrt{\lambda
_{v}\!-\!\phi
_{v}^{2}}\ln r\!\right)\!\right](\ln r)^{l}  \notag \\[.3pc]
&\quad\, +\sum_{v\in I_{3}}\sum_{l=0}^{2k_{v}-1}E_{l}~r^{-\phi
_{v}}(\ln r)^{l}
\end{align}
is the $r^{m}$ type  solution of the iterated equation~$(3)$.
Here{\rm ,} $A_{l},B_{l},$ $C_{l},D_{l},E_{l}$ are arbitrary
constants{\rm ,} $\phi_{v}$ is as given in $(6)$ and we divide the
index set $I=\{v=1,2,\dots, q\}$ into three parts{\rm :}
\begin{align*}
I_{1} &= \{ v\in I;\ \phi _{v}^{2}-\lambda _{v}>0\},\\[.3pc]
I_{2} &= \{ v\in I;\ \phi _{v}^{2}-\lambda _{v}<0\},\\[.3pc]
I_{3} &= \{ v\in I;\ \phi _{v}^{2}-\lambda _{v}=0\}.
\end{align*}
\end{theor}

\begin{proof} For any $v\in I,$ we can rewrite (7) as
\begin{equation*}
\left( \prod\limits_{j=1}^{q}L_{j}^{k_{j}}\right) (r^{m})=\left(
\beta _{v}^{k_{v}}(m)\prod\limits_{\substack{ j=1 \\ \text{ \
}j\neq v}}^{q}\beta_{j}^{k_{j}}(m)\right) r^{m}
\end{equation*}
or simply
\begin{equation}
\left( \prod\limits_{j=1}^{q}L_{j}^{k_{j}}\right) (r^{m})=\beta
_{v}^{k_{v}}(m)F(m).
\end{equation}
Here, we let $F(m)=\left( \prod\limits_{\substack{ j=1 \\ \text{\
}j\neq v}}^{q}\beta _{j}^{k_{j}}(m)\right) r^{m}$. Now, since
\begin{equation*}
\frac{\partial }{\partial m}\left(
\prod\limits_{j=1}^{q}L_{j}^{k_{j}}\right) (r^{m})=\left(
\prod\limits_{j=1}^{q}L_{j}^{k_{j}}\right) \left(\frac{\partial
}{\partial m}r^{m}\right)=\left(
\prod\limits_{j=1}^{q}L_{j}^{k_{j}}\right) (r^{m}\ln r),
\end{equation*}%
by taking the derivative with respect to $m$ on both sides of (9),
we get
\begin{align*}
\left( \prod\limits_{j=1}^{q}L_{j}^{k_{j}}\right) (r^{m}\ln r) &=
\frac{\partial }{\partial m}(\beta _{v}^{k_{v}}(m)F(m))\\[.3pc]
&= \beta _{v}^{k_{v}-1}(m) \{ k_{v}\beta
_{v}^{{\prime}}(m)F(m)+\beta _{v}(m)F^{{\prime }}(m)\}
\end{align*}
or simply
\begin{equation}
\left( \prod\limits_{j=1}^{q}L_{j}^{k_{j}}\right) (r^{m}\ln
r)=\beta_{v}^{k_{v}-1}(m)\Theta _{1}(m).
\end{equation}
Here, we set
\begin{equation*}
\Theta _{1}(m)=k_{v}\beta _{v}^{{\prime }}(m)F(m)+\beta
_{v}(m)F^{{\prime}}(m).
\end{equation*}
\end{proof}

Now, by taking the derivative with respect to $m$ on both sides of
(10), we obtain
\begin{equation*}
\left( \prod\limits_{j=1}^{q}L_{j}^{k_{j}}\right) (r^{m}(\ln
r)^{2})=\beta _{v}^{k_{v}-2}(m)\Theta _{2}(m),
\end{equation*}
where
\begin{equation*}
\Theta _{2}(m)=(k_{v}-1)\beta _{v}^{{\prime }}(m)\Theta
_{1}(m)+\beta_{v}(m)\Theta _{1}^{{\prime }}(m).
\end{equation*}

In a similar fashion, taking the successive derivatives $k_{v}-1$
times, with respect to $m$ on both sides of (9), we finally obtain
\begin{equation}
\left( \prod\limits_{j=1}^{q}L_{j}^{k_{j}}\right) (r^{m}(\ln
r)^{k_{v}-1})=\beta _{v}(m)\Theta _{k_{v}-1}(m).
\end{equation}
Here,
\begin{equation*}
\Theta _{k_{v}-1}(m)=2\beta _{v}^{{\prime }}(m)\Theta
_{k_{v}-2}(m)+\beta _{v}(m)\Theta _{k_{v}-2}^{{\prime }}(m).
\end{equation*}

Since the roots of the equation
\begin{equation*}
\beta _{v}(m)=m(m+2\phi _{v})+\lambda _{v}=0
\end{equation*}
are
\begin{equation*}
m_{v}^{(1)}=-\phi _{v}+\sqrt{\phi _{v}^{2}-\lambda _{v}}
\end{equation*}
and
\begin{equation*}
m_{v}^{(2)}=-\phi _{v}-\sqrt{\phi _{v}^{2}-\lambda _{v}},
\end{equation*}
we conclude from (11) that the functions
\begin{equation*}
r^{m_{v}^{(i)}}(\ln r)^{l}\qquad (i=1,2;\quad l=0,1,\dots,k_{v}-1)
\end{equation*}
are all solutions of eq.~(3). Thus, since the equation is linear,
by the superposition principle, the function
\begin{equation}
\sum_{v=1}^{q}\sum_{l=0}^{k_{v}-1}\left\{
A_{l}~r^{m_{v}^{(1)}}+B_{l}~r^{m_{v}^{(2)}}\right\} (\ln r)^{l}
\end{equation}%
is also a solution of (3).

We have three cases for the roots:

\begin{case}{\rm
If $v\in I_{1}$, then $m_{v}^{(1)}$ and $m_{v}^{(2)}$ are both
real. In this case, from (12), the function
\begin{equation}
\sum_{v\in I_{1}}^{{}}\sum_{l=0}^{k_{v}-1}r^{-\phi _{v}}\left\{
A_{l}~r^{\sqrt{\phi _{v}^{2}-\lambda _{v}}}+B_{l}~r^{-\sqrt{\phi
_{v}^{2}-\lambda _{v}}}\right\} (\ln r)^{l}.
\end{equation}
is a real-valued solution of (3).}
\end{case}

\begin{case}{\rm
If $v\in I_{2}$, then $m_{v}^{(1)}$ and $m_{v}^{(2)}$ are both
complex and conjugate. In this case, from (12), the function
\begin{equation}
\sum_{v\in I_{2}}^{{}}\sum_{l=0}^{k_{v}-1}r^{-\phi _{v}}\left\{
C_{l}~\cos \left(\sqrt{\lambda _{v}-\phi _{v}^{2}}\ln
r\right)+D_{l}\sin \left(\sqrt{\lambda _{v}-\phi _{v}^{2}}\ln
r\right)\right\} (\ln r)^{l}
\end{equation}
satisfies (3). Here, we use the Euler formula
\begin{align*}
r^{\pm i\sqrt{\lambda _{v}-\phi _{v}^{2}}} &= {\rm e}^{\pm
i\sqrt{\lambda _{v}-\phi _{v}^{2}}\ln r}\\[.3pc]
&= \left[ \cos \left(\sqrt{\lambda _{v}-\phi _{v}^{2}}\ln
r\right)\pm i\sin \left(\sqrt{\lambda _{v}-\phi _{v}^{2}}\ln
r\right)\right],
\end{align*}
and $C_{l}=A_{l}+B_{l},$ $D_{l}=i(A_{l}-B_{l})~$and $i=\sqrt{-1}$
as usual.}
\end{case}

\begin{case}{\rm
Finally, if $v\in I_{3}$, then $\ m_{v}^{(1)}=
m_{v}^{(2)}=-\phi_{v}$ is a multiple root. Thus, from (12), the
function
\begin{equation*}
\sum_{v\in I_{3}}^{{}}\sum_{l=0}^{k_{v}-1}\{
E_{l}~r^{m_{v}^{(1)}}\} (\ln r)^{l}
\end{equation*}
is a solution of (3). Now, from (9), since we have
\begin{equation*}
\left( \prod\limits_{j=1}^{q}L_{j}^{k_{j}}\right) (r^{m})=\beta
_{v}^{k_{v}}(m)F(m)=(m-m_{v}^{(1)})^{2k_{v}}F(m),
\end{equation*}\pagebreak

\noindent by taking the derivatives $2k_{v}-1$ times, with respect
to $m,$ on both sides of the above equality and letting
$m=m_{v}^{(1)},$ we obtain
\begin{equation*}
\left( \prod\limits_{v\in I_{3}}^{{}}L_{v}^{k_{v}}\right) (
r^{m_{v}^{(1)}}(\ln (r))^{l}) =0,\quad l=0,1,\dots,2k_{v}-1.
\end{equation*}
Hence, we conclude that the function
\begin{equation}
\sum_{v\in I_{3}}^{{}}\sum_{l=0}^{2k_{v}-1}E_{l}~r^{-\phi
_{v}}(\ln r)^{l}
\end{equation}
satisfies (3).

Summing up the above three cases with the superposition principle
we get (8), which proves the theorem.}\vspace{-.5pc}
\end{case}

\section{General solution for the iterated Euler equations}

In this section, we state the general solution of the iterated
Euler equations. In \cite{4}, for the Euler equation
\begin{equation*}
Eu=x^{2}\frac{{\rm d}^{2}u}{{\rm d}x^{2}}+\alpha x\frac{{\rm
d}u}{{\rm d}x}+\lambda u=0,
\end{equation*}
the general solutions for the iterated equations $E^{k}u=0$ are
given for any integer $k,$ where $\alpha$ and $\lambda$ are
arbitrary constants. Now consider the Euler equations
\begin{equation*}
E_{v}u=x^{2}\frac{{\rm d}^{2}u}{{\rm
d}x^{2}}+\alpha_{v}x\frac{{\rm d}u}{{\rm d}x}+\lambda _{v}u=0,
\end{equation*}
where $\alpha_{v}$ and $\lambda _{v}$ $(v=1,2,\dots,q)$ are
arbitrary constants.

The following result gives the general solutions of the iterated
Euler equations.

\begin{theor}[\!] The general solution of the iterated
Euler equations
\begin{equation*}
\left( \prod\limits_{v=1}^{q}E_{v}^{k_{v}}\right) u=(
E_{1}^{k_{1}}E_{2}^{k_{2}}\dots E_{q}^{k_{q}}) u=0
\end{equation*}
is
\begin{align*}
u &= \sum_{v\in I_{1}}\sum_{l=0}^{k_{v}-1}x^{-\phi _{v}}\left[
A_{l}~x^{\sqrt{\phi _{v}^{2}-\lambda _{v}}}+B_{l}~r^{-\sqrt{\phi
_{v}^{2}-\lambda _{v}}}\right] (\ln x)^{l} \\[.3pc]
&\quad\, +\!\sum_{v\in I_{2}}\sum_{l=0}^{k_{v}-1}x^{-\phi
_{v}}\!\left[\!C_{l}\cos \left(\!\!\sqrt{\lambda _{v}\!-\!\phi
_{v}^{2}}\ln
x\!\right)\!+\!D_{l}\sin \left(\!\!\sqrt{\lambda _{v}\!-\!\phi_{v}^{2}}\ln x\!\right)\!\right](\ln x)^{l} \\[.3pc]
&\quad\, +\sum_{v\in I_{3}}\sum_{l=0}^{2k_{v}-1}E_{l}~x^{-\phi
_{v}}(\ln x)^{l}.
\end{align*}
\end{theor}

\begin{proof}
In Theorem~1, by letting $n=1,$ and hence letting $\ r=x_{1}=x,$
$\alpha _{1}^{(v)}=\alpha _{v},$ we obtain the result for $\phi
_{v}=\frac{1}{2}(-1+\alpha _{v}).$\vspace{-.5pc}
\end{proof}


\begin{thebibliography}{9}
\bibitem{1} Alt\i n A, Solutions of type $r^{m}$ for a
class of singular equations, {\it Internat. J.~Math. and Math.
Sci.} {\bf 5(3)} (1982) 613--619

\bibitem{2} Alt\i n A, Some expansion formulas for a
class of singular partial differential equations, {\it Proc. Am.
Math. Soc.} {\bf 85(1)} (1982) 42--46

\bibitem{6} Lyakhov L N and Ryzhkov A V, Solutions of the $B$-polyharmonic
equation, {\it Differential Equations} {\bf 36(10)} (2000)
1507--1511, translated from {\it Differetsial'nye Uravneniya} {\bf
36(10)} (2000) 1365--1368

\bibitem{3} \"{O}zalp N and \c{C}etinkaya A, Radial
solutions of a class of iterated partial differential equations,
{\it Czechoslovak Math. J.} {\bf 55(130)} (2005) 531--541

\bibitem{4} \"{O}zalp N, $r^{m}$-type solutions
for a class of partial differential equations, {\it Commun. Fac.
Sci. Univ. Ank. Series A1} {\bf 49} (2001) 95--100

\bibitem{5} \"{O}zalp N and \c{C}etinkaya A, Expansion formulas
and Kelvin principle for a class of partial differential
equations, {\it Math. Balkanica (NS)} {\bf 15(3--4)} (2001)
219--226
\end{thebibliography}
\end{document}